\documentclass{ifacconf}

\usepackage{graphicx}      
\usepackage{natbib}        
\usepackage[utf8]{inputenc}
\usepackage{amsmath,bm}
\usepackage{amsfonts}
\usepackage{amssymb}
\usepackage{amsfonts}
\usepackage{caption}
\usepackage{subcaption}
\usepackage{algorithm}
\usepackage[noend]{algpseudocode}
\usepackage{xcolor}
\usepackage[normalem]{ulem}

\newtheorem{Theorem}{Theorem}
\numberwithin{Theorem}{section}

\newtheorem{Lemma}{Lemma}
\numberwithin{Lemma}{section}
\newtheorem{Problem}{Problem}
\numberwithin{Problem}{section}
\newtheorem{Corollary}{Corollary}
\newtheorem{Remark}{Remark}
\newtheorem{Assumption}{Assumption}
\usepackage{tikz}
	\usetikzlibrary{shapes,arrows}
	\tikzstyle{frame} = [draw, -latex]
	\tikzstyle{line} = [draw]
	\tikzstyle{line2} = [draw, dashdotted]
	\tikzstyle{line3} = [draw, dashed]
	\tikzstyle{line3UD} = [draw, dashed]
	\tikzstyle{place} = [circle, draw=black, fill=white, thick, inner sep=2pt, minimum size=1mm]
	\tikzstyle{place2} = [circle, draw=black, fill=black, thick, inner sep=2pt, minimum size=1mm]
	\tikzstyle{placeRed} = [circle, draw=red, fill=red, thick, inner sep=2pt, minimum size=1mm]
	\tikzstyle{vertex} = [circle, draw=black, fill=black, thick, inner sep=2pt, minimum size=1mm]
\usepackage{tkz-euclide}

\newcommand{\m}[1]{\mathbf{#1}}
\newcommand{\mc}[1]{\mathcal{#1}}
\newcommand{\mb}[1]{\mathbb{#1}}
\newcommand{\abs}[1]{\lVert{#1} \rVert}
\newcommand\numberthis{\addtocounter{equation}{1}\tag{\theequation}}
\newcommand*{\eprf}{\hfill\ensuremath{\blacksquare}}
\makeatletter
\def\algbackskip{\hskip-\ALG@thistlm}
\def\endthebibliography{%
  \def\@noitemerr{\@latex@warning{Empty `thebibliography' environment}}%
  \endlist
}
\makeatother
\begin{document}
\begin{frontmatter}

\title{Multi-agent Localization of A Common Reference Coordinate Frame: An Extrinsic Approach} 

\thanks[footnoteinfo]{The work of this paper has been supported by the National Research Foundation (NRF) of Korea under the grant NRF-2017R1A2B3007034.}

\author[First]{Quoc Van Tran} 
\author[First,Second]{Hyo-Sung Ahn} 

\address[First]{School of Mechanical Engineering,
        Gwangju Institute of Science and Technology, Gwangju, Republic of Korea (e-mail: $\{$tranvanquoc, hyosung$\}$@gist.ac.kr).}
\address[Second]{Department of Electrical Engineering, Colorado School of Mines, Golden, CO, USA (e-mail: hahn@mines.edu.}

\begin{abstract}                
This paper studies the problem of multi-agent cooperative localization of a common reference coordinate frame in $\mb{R}^3$. Each agent in a system maintains a body-fixed coordinate frame and its actual \textit{frame transformation} (translation and rotation) from the global coordinate system is unknown. The mobile agents aim to determine their \textit{trajectories of rigid-body motions} (or the frame transformations, i.e., rotations and translations) with respect to the global coordinate frame up to a common frame transformation by using local measurements and information exchanged with neighbors. We present two frame localization schemes which compute the rigid-body motions of the agents with asymptotic stability and finite-time stability properties, respectively. Under both localization laws, the estimates of the frame transformations of the agents converge to the actual frame transformations almost globally and up to an unknown constant transformation bias. Finally, simulation results are provided.
\end{abstract}


\end{frontmatter}

\section{Introduction}
\begin{figure}[t]
\centering
\begin{tikzpicture}
\node (py) at (0,1.,0) [label=below left:$\Sigma$]{};
\node (pz) at (0,0,1.5) [label=left:$$]{};
\node (px) at (1.,-0.2) [label=below:$$]{};
\node[place,scale=0.4] (pi) at (1.7,1.5){};
\node[place,scale=0.4] (pj) at (3.2,.3){};
\node[place,scale=0.01] (pi_y) at (1.2,2.4)[label=right:$^i\Sigma$]{};

\draw[black,rotate around={30:(pi)}] (0.7+0.5,0.9+0.3) rectangle (1.7+0.5,1.5+0.3);
\filldraw[black,rotate around={30:(0.9+0.6,1.5+.2)}] (.8+0.6,1.4+0.2) rectangle (1.1+0.6,1.6+0.2);
\filldraw[black,rotate around={30:(0.9+0.9,1.5-.3)}] (.8+0.9,1.4-.3) rectangle (1.1+0.9,1.6-.3);
\draw[{line width=1pt}] (0,0) [frame] -- (px);
\draw[{line width=1pt}] (0,0) [frame] -- (py);
\draw[{line width=1pt}] (0,0) [frame] -- (pz);

\draw[{line width=.5pt}] (pi) [frame] -- (2.5,1.94);
\draw[{line width=.5pt}] (pi) [frame] -- (pi_y);
\draw[{line width=0.5pt}] (pi) [frame] -- (.7,1.6);

\draw[{line width=1pt}, blue] (0,0) [frame] -- (pi)  node [pos =0.4, yshift=1ex, above] {$\m{T}_i$};
\draw[{line width=1pt}, blue] (0,0) [frame] -- (pj)  node [pos =0.5, yshift=-0ex, below] {$\m{T}_j$};
\draw[{line width=1pt}, blue] (pi) [frame] -- (pj)  node [pos =0.6, yshift=-0ex, above] {$\m{T}_{ij}$};

\draw[black,rotate around={-10:(pi)}] (0.82+2,0.9+-0.65) rectangle (1.82+2,1.5-0.65);
\filldraw[black,rotate around={-10:(0.95+2.2,1.6-1)}] (.85+2.2,1.5-1) rectangle (1.15+2.2,1.7-1);
\filldraw[black,rotate around={-10:(0.9+2.2,1.5-1.5)}] (.8+2.2,1.4-1.5) rectangle (1.1+2.2,1.6-1.5);
\draw[{line width=.5pt}] (pj) [frame] -- (3.2+0.1,.3+0.8) node [right] {$^j\Sigma$};
\draw[{line width=0.5pt}] (pj) [frame] -- (3.2+1,0.1);
\draw[{line width=0.5pt}] (pj) [frame] -- (3.2-0.6,-0.1);

\end{tikzpicture}
\caption{Relative transformation between two mobile agents $i$ and $j$ in $\mb{R}^3$. Each agent $i$ in the system measures relative frame transformations $\m{T}_{ij}\in SE(3)$ to its neighbors $j$s. The agents collaboratively estimate their actual rigid-body trajectories (or frame transformations) $\m{T}_i(t)\in SE(3)$.}
\label{fig:body_frame}
\end{figure}
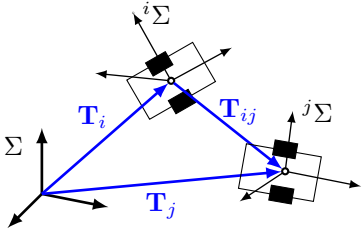

Consider a network of $n$ autonomous agents in $3$-dimensional space. Associated with each agent $i$,  there are a position vector, i.e., $\m{p}_i\in \mb{R}^3$, (normally correpsonding to its centroid) and a matrix in $SO(3)$, i.e., $\m{R}_i\in SO(3)$, its orientation, representing the orientation of its body-fixed coordinate frame, $^i\Sigma$, relative to the global coordinate frame $\Sigma$ (See also Fig. \ref{fig:body_frame}). This paper addresses the problem of estimating the trajectories of the rigid-body motions (as elements of the group of \textit{Euclidean transformations} $SE(3)$) or the time-varying poses of the agents which are characterized by the pairs $(\m{R}_i,\m{p}_i)\in SO(3)\times \mb{R}^3$. The global coordinate frame is unknown to all agents and they have only relative measurements and information communicated from their neighboring agents. To solve the problem, the agents in the system cooperatively localize a common reference coordinate frame and estimate their frame transformations with regard to the common frame.

In the literature, there has been recently a large number of works on consensus \cite{Igarashi2009, Sarltte2009, Thunberg2017,Markdahl2018tac, Gui2018auto, Zong2016, JWei2018tac} and estimation on $SO(3)$ \cite{Tron2014tac,Tron2016csm, Lee2017, Quoc2018cdc, Quoc2018tcns}. The consensus protocols have a wide range of applications in the orientation synchronization of rigid bodies in the Cartesian space. While the consensus algorithms designed directly on $SO(3)$ guarantee only local and asymptotic convergence property \cite{Igarashi2009, Sarltte2009}, those using local representations of orientations can provide almost global consensus and finite-time stability property \cite{Gui2018auto, Zong2016, JWei2018tac}. However, the local representations of orientations suffer from singularities, e.g., the angle-axis representation or the modified Rodriguez parameters, or the ambiguity in the orientation representations, e.g., the unit quaternions. For this reason, we only focus on the orientation control and estimation protocols which use directly the $SO(3)$ group to represent orientations. Orientation estimation approaches have been proposed recently and are widely used in network localization \cite{Tron2014tac,Tron2016csm, Lee2017} and formation control \cite{Lee2017, Quoc2018cdc, Quoc2018tcns}. The orientation estimation algorithms can guarantee almost global convergence of the computed orientations. 

The above-mentioned orientation estimation and control schemes can also be classified into intrinsic algorithms \cite{Igarashi2009, Tron2014tac, Markdahl2018tac} and extrinsic algorithms \cite{Sarltte2009, Lee2017, Thunberg2017, Quoc2018cdc, Quoc2018tcns}. In particular, the intrinsic algorithms design the orientation estimation and consensus laws directly on the Riemannian manifold, i.e., the sphere \cite{Markdahl2018tac} or the special orthogonal groups \cite{Igarashi2009,Tron2014tac}. By reshaping the cost function used in the estimation protocol, the convexity of the problem is guaranteed and the orientations of the agents (whose interaction graph is connected and undirected) can be estimated almost globally \cite{Tron2014tac, Markdahl2018tac}. Whereas, in the extrinsic approaches, rotation matrices are embedded into auxiliary matrices which are defined and evolve in the Euclidean ambient space through a typical consensus protocol. The auxiliary matrices are then exploited to control \cite{Sarltte2009} or estimate \cite{Lee2017, Thunberg2017, Quoc2018cdc, Quoc2018tcns} the orientation matrices. In contrast to the intrinsic algorithms, the extrinsic algorithms can guarantee almost global convergence of the rotation matrices for systems with general graph topologies which contain a spanning tree \cite{Lee2017,Quoc2018cdc, Quoc2018tcns}.

The consensus problem on $SE(3)$ was studied in \cite{Thunberg2016} with only local region of attraction. As the first contribution of this work, we formulate the \textit{frame localization problem} and propose an extrinsic-based algorithm to estimate the trajectories of the rigid-body motions of the agents in the system. In particular, the estimation law is designed by implementing consensus protocol on the auxiliary matrices in $\mb{R}^{4\times 4}$ and derived the poses of the agents from the auxiliary matrices. We show that the poses of the agents can be estimated almost globally and exponential fast up to an unknown constant frame transformation under the assumption of the existence of a spanning tree in the interaction graph. Secondly, a finite-time frame localization law is then proposed for systems with undirected and connected graph topologies. We establish almost global stability and the finite-time convergence of the estimated frame transformations of the agents to the actual frame transformations up to a common unknown transformation. The proposed frame localization protocols (on $SE(3)$) in this work are extended from the orientation estimation laws proposed in our previous works \cite{Quoc2018cdc, Quoc2018tcns}. Finally, simulation results are provided to show the effectiveness of the proposed frame localization schemes. 

The rest of this paper is outlined as follows. Preliminaries on the \textit{special Euclidean groups} and graph theory, and the problem formulation are presented in Section \ref{sec:preliminary}. Section \ref{sec:orient_local} proposes a frame localization law and establishes almost global exponential convergence of the estimated frame transformations. A frame localization scheme with finite-time stability property is introduced in Section \ref{sec:finite_time_frame_local}. We provide simulation results in Section \ref{sec:simulation}. Finally, Section \ref{sec:Conclusion} concludes this paper. 
\section{Preliminaries and Problem Formulation}\label{sec:preliminary} 
In this paper we use the following notations. 
Given two vectors $\m{x},\m{y}\in\mb{R}^3$, their dot product is denoted by $\m{x}\cdot\m{y}$. The symbol $\Sigma$ represents a global coordinate frame and the symbol $^k\Sigma$ with the superscript index $k$ denotes the $k$-th local coordinate frame. Let $\m{1}_n=[1,\ldots,1]^\top\in \mb{R}^n$ be the vector of all ones, and $\m{I}_3$ denotes the $3\times 3$ identity matrix. For two matrices $\m{A}$ and $\m{B}$, $\m{A}\otimes \m{B}$ denotes the Kronecker product between $\m{A}$ and $\m{B}$. The trace of a matrix is denoted by $\text{tr}(\cdot)$. For $\m{A}, \m{B}\in \mb{R}^{d\times l}$ the Frobenius metric is given by
$\abs{\m{A}-\m{B}}_F=\sqrt{\text{tr}\left\{(\m{A}-\m{B})^\top(\m{A}-\m{B})\right\}},$
which is the Euclidean distance in $\mb{R}^{d\times l}$.
\subsection{Special Euclidean Groups}
The set of rotation matrices in $\mb{R}^3$ is denoted by $SO(3)=\{\m{Q}\in\mb{R}^{3\times 3}~|~\m{Q}\m{Q}^\top=\m{I}_3,\text{det}(\m{Q})=1\}$. 
The space of $3\times 3$ skew-symmetric matrices is denoted by $\mathfrak{so}(3):=\{\m{A}\in \mb{R}^{3\times 3}|\m{A}^\top=-\m{A}\}$. For any $\boldsymbol{\omega}=[\omega_x,\omega_y,\omega_z]^\top\in \mb{R}^3$, the \textit{hat} map $(\cdot)^{\wedge}:~\mb{R}^3\rightarrow \mathfrak{so}(3)$ is defined such that $\boldsymbol{\omega}\times \m{x}=\boldsymbol{\omega}^\wedge\m{x},\forall \m{x}\in \mb{R}^3$, where
\begin{equation*}
\boldsymbol{\omega}^\wedge=\begin{bmatrix}
0 &-\omega_z &\omega_y\\
\omega_z &0 &-\omega_x\\
-\omega_y &\omega_x &0
\end{bmatrix}.
\end{equation*} 
The \textit{vee} map is the inverse of the \textit{hat} map and defined as $(\cdot)^\vee:~\mathfrak{so}(3)\rightarrow \mb{R}^3$ \cite{Bullo2005spr}. 

The \textit{special Euclidean group}, representing a trajectory of the motion (or poses, i.e., translation and rotation) of a rigid body agent, is given by a set of transformation matrices:
\begin{equation*}
SE(3)=\left\{ \m{T}=\begin{bmatrix}
\m{Q} &\m{p}\\
\m{0} &1
\end{bmatrix}\in \mb{R}^{4\times 4}\bigg| \m{Q}\in SO(3), \m{p}\in \mb{R}^3
\right\}.
\end{equation*} 
Note that the sets $SO(3)$ and $SE(3)$ are not vectorspaces, but they are \textit{matrix Lie groups} \cite{Barfoot2018}.
Let the set 
$$\mathfrak{se}(3):=\left\{ \boldsymbol{\xi}^\wedge=\begin{bmatrix}
\boldsymbol{\omega}^\wedge &\m{v}\\
\m{0} &0
\end{bmatrix}\in \mb{R}^{4\times 4}\bigg|~\omega^\wedge\in \mathfrak{so}(3), \m{v}\in \mb{R}^3\right\}.$$ The the \textit{hat} map $(\cdot)^{\wedge}:~\mb{R}^6\rightarrow \mathfrak{se}(3)$ and the \textit{vee} map $(\cdot)^\vee:~\mathfrak{se}(3)\rightarrow \mb{R}^6$ associated with $\mathfrak{se}(3)$ are given as
\begin{equation*}
\begin{bmatrix}
\m{v} \\ \boldsymbol{\omega}
\end{bmatrix}^\wedge :=
\begin{bmatrix}
\boldsymbol{\omega}^\wedge &\m{v}\\
\m{0} &0
\end{bmatrix} \text{ and }
\begin{bmatrix}
\boldsymbol{\omega}^\wedge &\m{v}\\
\m{0} &0
\end{bmatrix}^\vee:=\begin{bmatrix}
\m{v} \\ \boldsymbol{\omega}
\end{bmatrix}
,\end{equation*} respectively. The \textit{exponential map} relates the $\m{T}\in SE(3)$ and $\boldsymbol{\xi}^\wedge\in\mathfrak{se}(3)$ as
$$
\m{T}=\exp(\boldsymbol{\xi}^\wedge)=\sum_{k=0}^\infty\frac{1}{k!}(\boldsymbol{\xi}^\wedge)^k.
$$
The \textit{inverse} operator of the exponential map is given as
$$
\boldsymbol{\xi}=\text{ln}(\m{T})^\vee.
$$

\subsection{Graph theory}
An interaction graph characterizing an interaction topology of a multi-agent network is denoted by $\mc{G}=(\mc{V},\mc{E})$, where, $\mc{V}=\{1,\ldots,n\}$ denotes the vertex set and $\mc{E}\subseteq\mc{V}\times \mc{V}$ denotes the set of edges of $\mc{G}$. An edge is defined by the ordered pair $e_k=(i,j), k=1,\ldots,m,m=\vert \mathcal{E} \vert$. The graph $\mc{G}$ is said to be undirected if $(i,j)\in \mc{E}$ implies $(j,i)\in \mc{E}$, i.e. if $j$ is a neighbor of $i$, then $i$ is also a neighbor of $j$. If the graph $\mc{G}$ is directed, $(i,j)\in \mc{E}$ does not necessarily imply $(j,i)\in \mc{E}$. The set of neighboring agents of $i$ is denoted by $\mc{N}_i=\{j\in\mc{V}:(i,j)\in \mc{E}\}$. The Laplacian matrix $\m{L}=[l_{ij}]$ associated with $\mc{G}$ is defined as $l_{ij}=-1$ for $(i,j)\in \mc{E},~i\neq j$, $l_{ii}=-\sum_{j\in \mc{N}_i}l_{ij},~\forall i=1,\ldots,n$, and $l_{ij}=0$ otherwise. 
\subsection{Problem formulation}

Consider a network of $n$ mobile agents in $3$-dimensional space. 
Let $\m{p}_i$ and $\m{p}_i^i \in \mb{R}^3$ be the position of agent $i$ expressed in the global frame $\Sigma$ and its body-fixed coordinate frame $^i\Sigma$, respectively. 
The pair $(\m{R}_i,\m{p}_i) \in SO(3)\times \mb{R}^3$ characterizes the \textit{pose} of each agent $i$ in the Cartesian ambient space. The rigid body motion of agent $i$ (or the $i$-th frame transformation) is given as $\m{T}_i(t)=\begin{bmatrix}
\m{R}_i(t) &\m{p}_i(t)\\
\m{0} &1
\end{bmatrix}\in SE(3)$ 
and its inverse transformation can be computed as
$$\m{T}_i^{-1}=\begin{bmatrix}
\m{R}_i^\top &-\m{R}_i^\top\m{p}_i\\
\m{0} &1
\end{bmatrix}\in SE(3).$$ The relative transformation between the two corresponding body-fixed coordinate frames of agents $i$ and $j$, which is denoted as $\m{T}_{ij}\in SE(3)$, is given as
\begin{equation}\label{eq:relative_transform}
\m{T}_{ij} = \m{T}_i^{-1}\m{T}_j=\begin{bmatrix}
\m{R}_i^\top\m{R}_j &\m{R}_i^\top(\m{p}_j-\m{p}_i)\\
\m{0} &1
\end{bmatrix}.
\end{equation}
Let $\m{R}_{ij}:=\m{R}_i^\top\m{R}_j$ be the relative orientation between two local coordinate frames $^i\Sigma$ and $^j\Sigma$ and $\m{p}_{ij}^i:=\m{R}_i^\top\m{p}_{ij}=\m{R}_i^\top(\m{p}_j-\m{p}_i)$ the relative position between agent $i$ and $j$ which are measured locally in the local frame of agent $i$. Then, the relative transformation can be expressed as $\m{T}_{ij}=\begin{bmatrix}
\m{R}_{ij} &\m{p}_{ij}^i\\
\m{0} &1
\end{bmatrix}$. 

The kinematic of the rigid body motion of agent $i$ is given as
\begin{equation}
\dot{\m{T}}_i=\m{T}_i\begin{bmatrix}
\m{v}^i_i \\ \boldsymbol{\omega}^i_i
\end{bmatrix}^\wedge =\begin{bmatrix}
\m{R}_i &\m{p}_i\\
\m{0} &1
\end{bmatrix}
\begin{bmatrix}
{\boldsymbol{\omega}^i_i}^\wedge &\m{v}^i_i\\
\m{0} &0
\end{bmatrix}=\begin{bmatrix}
\m{R}_i{\boldsymbol{\omega}^i_i}^\wedge &\m{R}_i\m{v}^i_i\\
\m{0} &0
\end{bmatrix}
\end{equation}
where $\m{v}^i_i\in \mb{R}^3$ and $\boldsymbol{\omega}^i_i\in \mb{R}^3$ denote the linear velocity and the angular velocity of agent $i$ measured in $^i\Sigma$. We assume that each agent $i$ is able to measure $\m{v}^i_i$ and $\boldsymbol{\omega}^i_i$ and the relative transformations (translation and rotation) to its neighboring agents without noise. If an edge $(i,j)\in \mc{E}$, then agent $i$ can measure $\m{T}_{ij}\in SE(3)$ and it also can receive information communicated from agent $j$. For this, we make the following assumptions.

\begin{Assumption}\label{ass:measurements}
Each agent $i$ in the system locally measures its body velocity, i.e., $[\m{v}_i^i,\boldsymbol{\omega}_i^i]\in \mb{R}^3\times \mb{R}^3$, and the relative transformation $\m{T}_{ij}\in SE(3)$ defined in \eqref{eq:relative_transform} with regard to its neighbors $j\in \mc{N}_i$.
\end{Assumption}
\begin{Assumption}\label{ass:graph_topology}
The underlying interaction graph $\mc{G}(\mc{V},\mc{E})$ contains a spanning tree.
\end{Assumption}
The agents in the system aim to estimate for their actual rigid body motions $\m{T}_i(t)\in SE(3),i=1,\ldots,n,$, a process is called \textit{frame localization}, and each agent $i$ holds an estimate of the body transformation $\hat{\m{T}}_i=\hat{\m{T}}_i(\hat{\m{R}}_i,\hat{\m{p}}_i)\in SE(3), \forall i\in \mc{V}$. By using the local measurements in Assumption \ref{ass:measurements}, the objective is for the agents to cooperatively localize the  global coordinate frame, e.g., $\Sigma$, up to a transformation, $\m{T}_c(\m{R}_c,\m{p}_c)\in SE(3)$, which is unknown but deterministic and common to all agents.
\begin{Problem}[Asymptotic Frame Localization]\label{prob:network_localization}
Consider a system of $n$ mobile agents in $\mb{R}^3$. Under the Assumptions \ref{ass:measurements} and \ref{ass:graph_topology}, design a cooperative localization scheme for each agent $i$ to estimate its transformation $\m{T}_i$ up to an unknown constant transformation $\m{T}_c$.
\end{Problem}
\begin{Problem}[Finite-Time Frame Localization]\label{prob:finite_time_network_localization}
Consider a system of $n$ mobile agents in $\mb{R}^3$. Under the Assumptions \ref{ass:measurements} and assume that $\mc{G}$ is connected and undirected, design a cooperative localization scheme for each agent $i$ to estimate its transformation $\m{T}_i$ up to an unknown constant transformation $\m{T}_c$ in a finite time.
\end{Problem}
\section{Distributed Esimation of A Common Reference Frame}\label{sec:orient_local}
This section presents a distributed estimation protocol and establishes the almost global asymptotic convergence of the estimated poses to the actual poses of the agents up to a common reference transformation by using the relative pose measurements.
\subsection{Propose Estimation Law}
For each agent $i$ we introduce an auxiliary matrix $\m{P}_i\in \mb{R}^{4\times 3}$ as follows
\begin{equation}\label{eq:P_i}
\m{P}_i := \begin{bmatrix}
\m{Q}_i &\m{q}_i\\
\m{0} &1
\end{bmatrix},
\end{equation}
where $\m{Q}_i\in \mb{R}^{3\times 3}$ is a nonsingular matrix and $\m{q}_i\in \mb{R}^3$. Note that $\m{P}_i$ is defined in the Cartesian ambient space and $\m{Q}_i(0)$ has full-column rank and initialized randomly. Note that the set of nonsingular matrices in $\mb{R}^{3\times 3}$ is a dense set of the set of $3\times 3$ matrices, i.e., if $\m{Q}_i$ is initialized randomly from a continuous uniform distribution on its entries, then $\m{Q}_i$ will be almost surely nonsingular. Each agent $i$ implements the following localization law
\begin{equation}\label{eq:asymptotic_law}
\dot{\m{P}}_i(t)=-\begin{bmatrix}
\m{v}^i_i \\ \boldsymbol{\omega}^i_i
\end{bmatrix}^\wedge\m{P}_i(t)+ \sum_{j\in \mc{N}_i}\Big(\m{T}_{ij}(t)\m{P}_j(t)-\m{P}_i(t)\Big)
\end{equation} 
where $\m{P}_j\in \mb{R}^{4\times 4}$ is the auxiliary matrix associated with agent $j$ and it is communicated from agent $j\in \mc{N}_i$. In contrast to the intrinsic algorithms in the literature \cite{}, the frame localization law \eqref{eq:asymptotic_law} evolves in the Cartesian ambient space and the frame transformation estimate of each agent, $\hat{\m{T}}_i$, is derived from the corresponding auxiliary matrix $\m{P}_i$. The way, which constructs $\hat{\m{T}}_i$ from $\m{P}_i$, will be described latter in this section.
\subsection{Analysis}
The localization law \eqref{eq:asymptotic_law} can be rewritten as
\begin{equation*}
\m{T}_i\dot{\m{P}}_i=-\m{T}_i\begin{bmatrix}
\m{v}^i_i \\ \boldsymbol{\omega}^i_i
\end{bmatrix}^\wedge\m{P}_i+ \sum_{j\in \mc{N}_i}\Big(\m{T}_{j}\m{P}_j-\m{T}_i\m{P}_i\Big).
\end{equation*} 
By introducing the transformation $\m{P}_i=\m{T}_i^{-1}\m{S}_i$ and noticing that $\dot{\m{S}}_i=\m{T}_i\dot{\m{P}}_i+\dot{\m{T}}_i\m{P}_i=\m{T}_i\dot{\m{P}}_i+\m{T}_i\begin{bmatrix}
\m{v}^i_i \\ \boldsymbol{\omega}^i_i
\end{bmatrix}^\wedge\m{P}_i$. Therefore, the above equation can be expressed as
\begin{equation}
\dot{\m{S}}_i=\sum_{j\in \mc{N}_i}(\m{S}_j-\m{S}_i).
\end{equation}
Let $\m{S}:=[\m{S}_1^\top,\ldots,\m{S}_n^\top]^\top\in \mb{R}^{4n\times 4}$ be the stack matrix of all $\m{S}_i\in \mb{R}^{4\times 4},i=1,\ldots,n$. By combining the above frame localization dynamics for all agents we obtain a compact form
\begin{equation}\label{eq:localization_compact_form}
\dot{\m{S}}=-(\m{L}\otimes \m{I}_4)\m{S}.
\end{equation}
\begin{Theorem}\label{thm:exp_converg_S_i}
Assume that Assumptions \ref{ass:measurements} and \ref{ass:graph_topology} hold. Under the frame localization law \eqref{eq:asymptotic_law}, $\m{S}(t)$ in \eqref{eq:localization_compact_form} globally exponentially converge to 
\begin{equation*}
(\m{1}_n\otimes \m{I}_4)(\m{w}_1\times \m{I}_4)\m{S}(0),
\end{equation*}
where $\m{w}_1=[w_{11},\ldots,w_{1n}],\m{w}_1\m{1}_n=1,$ is the left eigenvector of the Laplacian $\m{L}$ corresponding to the zero eigenvalue.
\end{Theorem}
\begin{pf}
Since $\mc{G}$ has a spanning tree the associated Laplacian $\m{L}\in \mb{R}^{n\times n}$ has a simple zero eigenvalue and the other eigenvalues have positive real parts. The right and left eigenvectors corresponding to the zero eigenvalue are $\m{1}_n$ and $\m{w}_1=[w_{11},\ldots,w_{1n}],w_{1i}>0,\m{w}_1\m{1}_n=1,$ respectively \cite[Lemma 1]{Ren2004}. Further, there exists $\m{S}_c\in \mb{R}^{4\times 4}$ such that $\m{S}(t)$ globally exponentially converges to $(\m{1}_n\otimes \m{I}_4)\m{S}_c$.

Consider the solution to \eqref{eq:localization_compact_form} as
\begin{equation*}
\m{S}(t)=\exp^{-(\m{L}\otimes \m{I}_4)t}\m{S}(0).
\end{equation*}
Let the Jordan form of $\m{L}$ be $\m{L}=\m{UDV}^{-1}$ where $\m{D}=\text{diag}\{0,\lambda_2,\ldots,\lambda_n\}$ whose diagonal terms are eigenvalues of $\m{L}$, $\m{U}=[\m{1}_n,\m{r}_2,\ldots,\m{r}_n]$ and $\m{V}^{-1}=[\m{w}_1^\top,\ldots,\m{w}_n^\top]^\top\in \mb{R}^{n\times n}$. Then the steady-state solution $\lim_{t\rightarrow\infty}\m{S}(t)=(\m{1}_n\otimes \m{I}_4)(\m{w}_1\times \m{I}_4)\m{S}(0)$. Thus, $\m{S}_i,i\in \mc{V},$ converge to $\m{S}_c:=(\m{w}_1\times \m{I}_4)\m{S}(0)$, i.e., a convex combination of the initial matrices $\{\m{S}_i(0)\}_{i\in \mc{V}}$. \eprf
\end{pf}
The steady-state matrix $\m{S}_c\in \mb{R}^{4\times 4}$ is given as
\begin{align*}
\m{S}_c&=(\m{w}_1\times \m{I}_4)\m{S}(0)=\sum_{i=1}^n{w}_{1i}\m{T}_i(0)\m{P}_i(0)\\
&=\sum_{i=1}^n{w}_{1i}\begin{bmatrix}
\m{R}_i(0)\m{Q}_i(0) &\m{R}_i(0)\m{q}_i(0)+\m{p}_i(0)\\
\m{0} &1
\end{bmatrix}\\
&=\begin{bmatrix}
\m{Q}_c &\m{q}_c\\
\m{0} &1
\end{bmatrix},\numberthis \label{eq:S_c}
\end{align*}
where $\m{Q}_c:=\sum_{i=1}^n{w}_{1i}\m{R}_i(0)\m{Q}_i(0)\in \mb{R}^{3\times 3}$, $\m{q}_c:=\sum_{i=1}^n{w}_{1i}\Big(\m{R}_i(0)\m{q}_i(0)+\m{p}_i(0)\Big)\in \mb{R}^3$, and ${w}_{1i}$ denotes the $i$-th entry of the left eigenvector $\m{w}_1$.

At a time instant $t$, let $\m{S}_i(t)=\begin{bmatrix}
\m{Q}_{S_i}(t)\in \mb{R}^{3\times 3} &\m{q}_{S_i}(t)\in \mb{R}^3\\
\m{0} &1
\end{bmatrix}$ then the auxiliary matrix $\m{P}_i$ is computed as
\begin{equation}\label{eq:relation_Pi_Si}
\m{P}_i=\m{T}_i^{-1}\m{S}_i=\begin{bmatrix}
\m{R}_i^\top\m{Q}_{S_i} &\m{R}_i^\top(\m{q}_{S_i}-\m{p}_i)\\
\m{0} &1
\end{bmatrix}.
\end{equation}
Since $\m{S}_i\rightarrow\m{S}_c,\forall i \in \mc{V},$ globally exponentially as $t\rightarrow\infty$ (Theorem \ref{thm:exp_converg_S_i}) we have the following lemma.
\begin{Lemma}\label{lm:steady_state_Pi}
Assume that Assumptions \ref{ass:measurements} and \ref{ass:graph_topology} hold. Under the frame localization law \eqref{eq:asymptotic_law}, $$\m{P}_i(t)\rightarrow\begin{bmatrix}
\m{R}_i^\top\m{Q}_c &\m{R}_i^\top(\m{q}_c-\m{p}_i)\\
\m{0} &1
\end{bmatrix},\forall i \in \mc{V},$$ globally exponentially as $t\rightarrow\infty$, i.e., $\m{Q}_i\rightarrow\m{R}_i^\top\m{Q}_c$ and $\m{q}_i\rightarrow\m{R}_i^\top(\m{q}_c-\m{p}_i)$, where $\m{Q}_c$ and $\m{q}_c$ are defined in \eqref{eq:S_c}.
\end{Lemma}
\subsection{Construction of Frame Transformations}
We now assume that estimates of orientation, $\hat{\m{R}}_i$, and position, $\hat{\m{p}}_i$, of agent $i$ are derived from $\m{Q}_i$ and $\m{p}_i$ as follows ($\m{Q}_i$ and $\m{p}_i$ are defined in \eqref{eq:P_i}). The orientation estimate $\hat{\m{R}}_i^T(t)$ is constructed from $\m{Q}_i(t)$ by the Gram-Schmidt procedure (GSOP, see Appendix \ref{app:gsop}) and 
\begin{equation}\label{eq:computated_pi_hat}
\hat{\m{p}}_i(t):=-\hat{\m{R}}_i(t)\m{q}_i(t).
\end{equation}
It is noticed from Lemma \ref{lm:steady_state_Pi} that 
\begin{equation}\label{eq:hat_pi_i}
\hat{\m{p}}_i^i(\infty):=-\m{q}_i(\infty)=\m{R}_i^\top(\m{p}_i-\m{q}_c)=\m{p}_i^i-\m{R}_i^\top\m{q}_c
\end{equation}
specifies the estimate of position of agent $i$ expressed locally in body frame $^i\Sigma$. It follows that the position of agent $i$ expressed in $^i\Sigma$, i.e., ${\m{p}}_i^i$, is estimated up to a common constant translation $\m{q}_c\in \mb{R}^3$.

Let $\m{Z}_0:=\left[\big(\m{R}_1(0)\m{Q}_1(0)\big)^\top,\ldots,\big(\m{R}_n(0)\m{Q}_n(0)\big)^\top\right]^\top\in \mb{R}^{3n\times 3}$. Then, we have the following result.
\begin{Corollary}\label{corol:almost_global_convergence}
Assume that Assumptions \ref{ass:measurements} and \ref{ass:graph_topology} hold. Under the frame localization law \eqref{eq:asymptotic_law}, if $\hat{\m{R}}_i^T(t)$ is constructed from $\m{Q}_i(t)$ by the Gram-Schmidt procedure (GSOP) and $\hat{\m{p}}_i$ is computed by \eqref{eq:computated_pi_hat}, then there exist an unknown constant transformation  
\begin{equation}\label{eq:Tc}
\m{T}_c := \begin{bmatrix}
\m{R}_c &\m{q}_c\\
\m{0} &1
\end{bmatrix}\in SE(3)
\end{equation} 
such that $\hat{\m{T}}_i(t)\rightarrow\m{T}_c^{-1}\m{T}_i(t)$ as $t\rightarrow\infty$, for all $i\in \mc{V}$, if $\text{Range}(\m{Z}_0)\cap\text{null}(\m{w}_1\otimes \m{I}_3)=\emptyset$.
\end{Corollary}
\begin{pf}
Let $\m{\hat{R}}_i^T$ and $\m{R}_{Si}$ be derived from $\m{Q}_{i}$ and $\m{Q}_{Si}$ by the GSOP, respectively. It follows from Lemma \ref{lm:invariantOrthogonalization} in Appendix that $\m{\hat{R}}_i^T=\m{R}_i^T\m{R}_{Si}$ for all $t$. Since $\m{Q}_i\rightarrow\m{R}_i^\top\m{Q}_c$ (Lemma \ref{lm:steady_state_Pi}), $\m{\hat{R}}_i^\top\rightarrow\m{R}_i^\top\m{R}_{c},\forall i\in\mc{V},$ as $t\rightarrow \infty$, where the unknown constant orientation $\m{R}_{c}:=\text{GSOP}(\m{Q}_c)$. As a result, from \eqref{eq:computated_pi_hat} and Lemma \ref{lm:steady_state_Pi}, we have
$$
\hat{\m{p}}_i\rightarrow -\hat{\m{R}}_i\m{R}_i^\top(\m{q}_c-\m{p}_i)=\m{R}_c^\top(\m{p}_i-\m{q}_c),\forall i\in \mc{V},
$$
as $t\rightarrow\infty$. Consequently, one has
$$
\lim_{t\rightarrow\infty}\hat{\m{T}}_i(t)=\begin{bmatrix}
\m{R}_{c}^\top\m{R}_i &\m{R}_c^\top(\m{p}_i-\m{q}_c)\\
\m{0} &1
\end{bmatrix}=\m{T}_c^{-1}\m{T}_i(t),
$$
where $\m{T}_c\in SE(3)$ is defined in \eqref{eq:Tc}.

For the validity of the estimated frame transformations, $\hat{\m{T}}_i(\m{\hat{R}}_i,\hat{\m{p}}_i)$, the singularity of the steady-state matrix $\m{Q}_c$ defined in \eqref{eq:S_c} is undesired. For this, we now show that $\m{Q}_c$ is nonsingular if the initial matrices satisfy the condition $\text{Range}(\m{Z}_0)\cap\text{null}(\m{w}_1\otimes \m{I}_3)=\emptyset$. From \eqref{eq:S_c}, $\m{Q}_c$ is explicitly computed as 
\begin{align*}
\m{Q}_c&=\sum_{i=1}^n{w}_{1i}\m{R}_i(0)\m{Q}_i(0)=(\m{w}_1\otimes \m{I}_3)\m{Z}_0\\
&=\Big[(\m{w}_1\otimes \m{I}_3)[\m{Z}_0]_1,(\m{w}_1\otimes \m{I}_3)[\m{Z}_0]_2,(\m{w}_1\otimes \m{I}_3)[\m{Z}_0]_3 \Big]
\end{align*}
where $[\m{Z}_0]_i$ denotes the $i$-th column of $\m{Z}_0$. It follows that $\m{Q}_c$ contains linearly independent columns if and only if $(\m{w}_1\otimes \m{I}_3)[\m{Z}_0]_i\neq\m{0},\forall i=1,2,3,$ and column vectors of $\m{Z}_0$ are linearly independent. The second condition follows from the nonsingularity of $\m{Q}_i(0)$ (for almost all random initializations of the entries of $\m{Q}_i(0)$) and the first condition implies $\text{Range}(\m{Z}_0)\cap\text{null}(\m{w}_1\otimes \m{I}_3)=\emptyset$. \eprf
\end{pf}
\begin{Corollary}
Since the dimension of $\text{null}(\m{w}_1\otimes \m{I}_3)$ is $3(n-1)$ which is a lower-dimensional subspace of $\mb{R}^{3n}$ and hence its Lebesgue measure is zero. Thus, the steady-state estimates of the frame transformations $\hat{\m{T}}_i,i\in \mc{V},$ are well-defined for almost all initial matrices $\{\m{P}_i(0)\}_{i\in \mc{V}}$. Further, the frame transformations of the agents are computed almost globally exponentially up to a common constant transformation $\m{T}_c(\m{R}_c,\m{q}_c)$. In other words, the frame transformations of the agents are computed relative to a reference frame whose frame transformation is $\m{T}_c(\m{R}_c,\m{q}_c)$.
\end{Corollary}
\begin{Remark}
Though the steady-state estimates of the frame transformations of the agents are proper they might not be well-defined at some time instants (See also \cite{Quoc2019tcns}). Indeed, if $\{\m{Q}_i\}_{i\in \mc{V}}$ in \eqref{eq:P_i} are initialized randomly in $\m{R}^{3\times 3}$ then some of $\{\m{R}_i(0)\m{Q}_i(0)\}_{i\in\mc{V}},$ have negative and some of those have positive determinants. Since all $\{\m{R}_i(t)\m{Q}_i(t)\}_{i\in\mc{V}},$ converge to $\m{Q}_c$, at least one of $\{\m{R}_i(t)\m{Q}_i(t)\}_{i\in\mc{V}},$ whose determinant changes sign. Thus, its determinant becomes zero at some time instants. Consequently, $\m{Q}_i$ is nonsingular at those points. 
\end{Remark}
\begin{Remark}
If we instead construct the first two columns of $\hat{\m{R}}_i^\top$ from the first two column vectors of $\m{Q}_i$ by the Gram-Schmidt orthonormalization process and the third column vector of $\hat{\m{R}}_i^\top$ is the cross product of the first two column vectors, then it can be shown that $\hat{\m{R}}_i$ is well-defined for all $t>0$ for almost all initial matrices $\{\m{Q}_i(0)\}_{i\in \mc{V}}$. Indeed, it is equivalent to show that the first two column vectors of $\m{S}_i(t)$ in \eqref{eq:localization_compact_form} are linearly independent for all initial matrices $\{\m{S}_i(0)\}_{i\in \mc{V}}$ but a set of zero measure \cite{Quoc2019tcns}.
\end{Remark}

The frame localization scheme with asymptotic convergence property is illustrated in Algorithm \ref{alg:asympt_estimator}
\begin{algorithm}[t]
\begin{algorithmic}[1]
\BState \emph{Initialization}: $t\gets0, ~\m{P}_i(\m{Q}_i(0),\m{q}_i(0))\in \mb{R}^{4\times 4}$ as \eqref{eq:P_i}.
\BState \emph{Estimation loop}:
\State $~$\textbf{repeat}
\State $\quad$\textbf{for all} $i\in\mc{V}$ \textbf{do}
\State $\qquad\m{P}_i(t) \gets$ integrate \eqref{eq:asymptotic_law} 
\State $\qquad\m{\hat{R}}^T_i \gets$ \text{GSOP}($\m{Q}_i$)
\State $\qquad\hat{\m{p}}_i(t) \gets$ $-\hat{\m{R}}_i(t)\m{q}_i(t)$
\State $\quad$\textbf{end for}
\State $~$\textbf{until finish.}
\BState \emph{End estimation loop}.
\end{algorithmic}
\caption{Almost Global Asymptotic Frame Localization}
\label{alg:asympt_estimator}
\end{algorithm}
\section{Finite-time Frame Localization}\label{sec:finite_time_frame_local}
In this section, a finite-time frame localization law is proposed for systems with undirected and connected graph topologies. We establish an almost global stability and the finite-time convergence of the estimated frame transformations of the agents to the actual frame transformations up to an unknown common transformation.
\subsection{Proposed Finite-time Frame Localization Law}
Each agent $i$ holds an auxiliary matrix $\m{P}_i\in \mb{R}^{4\times 4}$ as defined in \eqref{eq:P_i} and the estimate of the frame transformation of agent $i$ is constructed from $\m{P}_i$ by following the same computations in Corollary \ref{corol:almost_global_convergence}. For each agent $i$, we propose the following frame localization law
\begin{equation}\label{eq:finite-time_law}
\dot{\m{P}}_i(t)=-\begin{bmatrix}
\m{v}^i_i \\ \omega^i_i
\end{bmatrix}^\wedge\m{P}_i(t)+ \sum_{j\in \mc{N}_i}\frac{\m{T}_{ij}(t)\m{P}_j(t)-\m{P}_i(t)}{||\m{T}_{ij}(t)\m{P}_j(t)-\m{P}_i(t)||_{F}^\alpha},
\end{equation} 
where the positive scalar $0<\alpha<1$. The above frame localization law is continuous due to the Remark \ref{rmk:continuous_finite_time_law}.

To analyze the above localization law, the following Lemma is useful. 
\begin{Lemma}\label{eq:equiv_denominator}
The denominator of the second term in the right hand side of \eqref{eq:finite-time_law} can be equivalently computed as
\begin{equation}\label{eq:finite-time_rewriten_denominator}
||\m{T}_{ij}\m{P}_j-\m{P}_i||_{F}^\alpha=||\m{T}_{j}\m{P}_j-\m{T}_{i}\m{P}_i||_{F}^\alpha
\end{equation}
\end{Lemma}
\begin{pf}
First, we have
\begin{align*}
&\m{T}_{ij}\m{P}_j-\m{P}_i=\\
&\begin{bmatrix}
\m{R}_i^\top\m{R}_j &\m{R}_i^\top(\m{p}_j-\m{p}_i)\\
\m{0} &1
\end{bmatrix}\begin{bmatrix}
\m{Q}_j &\m{q}_j\\
\m{0} &1
\end{bmatrix}
-\begin{bmatrix}
\m{Q}_i &\m{q}_i\\
\m{0} &1
\end{bmatrix}\\
&=\begin{bmatrix}
\m{R}_i^\top\m{R}_j\m{Q}_j &\m{R}_i^\top\m{R}_j\m{q}_j+\m{R}_i^\top(\m{p}_j-\m{p}_i)\\
\m{0} &1
\end{bmatrix}-\m{P}_i\\
&=\begin{bmatrix}
\m{R}_i^\top &\m{0}\\
\m{0} &1
\end{bmatrix}\begin{bmatrix}
\m{R}_j\m{Q}_j &\m{R}_j\m{q}_j+(\m{p}_j-\m{p}_i)\\
\m{0} &1
\end{bmatrix}-\m{P}_i\\
&=\begin{bmatrix}
\m{R}_i^\top &\m{0}\\
\m{0} &1
\end{bmatrix}\begin{bmatrix}
\m{R}_j\m{Q}_j &\m{R}_j\m{q}_j+\m{p}_j\\
\m{0} &1
\end{bmatrix}-\begin{bmatrix}
\m{Q}_i &\m{q}_i+\m{R}_i^\top\m{p}_i\\
\m{0} &1
\end{bmatrix}\\
&=\begin{bmatrix}
\m{R}_i^\top &\m{0}\\
\m{0} &1
\end{bmatrix}\m{T}_j\m{P}_j-\begin{bmatrix}
\m{R}_i^\top &\m{0}\\
\m{0} &1
\end{bmatrix}\begin{bmatrix}
\m{R}_i\m{Q}_i &\m{R}_i\m{q}_i+\m{p}_i\\
\m{0} &1
\end{bmatrix}\\
&=\m{K}_i(\m{T}_j\m{P}_j-\m{T}_i\m{P}_i),
\end{align*}
where $\m{K}_i:=\begin{bmatrix}
\m{R}_i^\top &\m{0}\\
\m{0} &1
\end{bmatrix}$.
By using this relation, one has
\begin{align*}
&||\m{T}_{ij}\m{P}_j-\m{P}_i||_{F}^\alpha=\text{tr}\Big[(\m{T}_{ij}\m{P}_j-\m{P}_i)^\top(\m{T}_{ij}\m{P}_j-\m{P}_i) \Big]\\
&=\text{tr}\Big[(\m{T}_j\m{P}_j-\m{T}_i\m{P}_i)^\top\m{K}_i^\top\m{K}_i(\m{T}_j\m{P}_j-\m{T}_i\m{P}_i) \Big]\\
&=\text{tr}\Big[(\m{T}_j\m{P}_j-\m{T}_i\m{P}_i)^\top(\m{T}_j\m{P}_j-\m{T}_i\m{P}_i) \Big]\\
&=||\m{T}_{j}\m{P}_j-\m{T}_{i}\m{P}_i||_{F}^\alpha,
\end{align*}
which competes the proof. \eprf
\end{pf}
\subsection{Analysis}
By using the transformation $\m{P}_i=\m{T}_i^{-1}\m{S}_i,~\m{S}_i\in \mb{R}^{4\times 4},$ and the above Lemma, the frame localization law \eqref{eq:finite-time_law} can be written as
\begin{equation}\label{eq:finite_transformed_system}
\dot{\m{S}}_i=\sum_{j\in \mc{N}_i}\frac{\m{S}_j-\m{S}_i}{||\m{S}_j-\m{S}_i||_F^\alpha}.
\end{equation}
\begin{Remark}\label{rmk:continuous_finite_time_law}
The above system is time-continuous as will be shown in the following. Let $[\m{S}_i]_k\in \mb{R}^4$ be the $k$-th column vector of $\m{S}_i$ and $\m{s}_i:=[[\m{S}_i]_1^\top,[\m{S}_i]_2^\top,\ldots,[\m{S}_i]_4^\top]\in \mb{R}^{16}$ be the stacked vector of all column vectors of $\m{S}_i$, for all $i\in \{1,\ldots,n\}$. Furthermore, from the definition of the Frobenius norm, we have
\begin{align*}
\abs{\m{S}_i-\m{S}_j}_F&=\sqrt{\text{tr}\{(\m{S}_i-\m{S}_j)^\top(\m{S}_i-\m{S}_j)\}}\\
&=\sqrt{(\m{s}_i-\m{s}_j)^\top(\m{s}_i-\m{s}_j)}\\
&= \abs{\m{s}_i-\m{s}_j},
\end{align*} 
where $\abs{\cdot}$ denotes the Euclidean norm. By using the above equation, we can rewrite \eqref{eq:finite_transformed_system} into a vector form as
\begin{equation}
\dot{\m{s}}_i(t)=\sum_{j\in\mc{N}_i}\frac{\m{s}_j-\m{s}_i}{\abs{\m{s}_i-\m{s}_j}^\alpha},
\end{equation}
with $\alpha\in (0,1)$ the right hand side of the above equation is continuous \cite{Minh2017cdc}. If $\alpha \geq 1$, it is discontinuous \cite{Cortes2006}. \eprf
\end{Remark}

Let $\m{S}:=[\m{S}_1^\top,\ldots,\m{S}_n^\top]^\top\in \mb{R}^{4n\times 4}$ be the stack matrix of all $\m{S}_i\in \mb{R}^{4\times 4},i=1,\ldots,n$. By combining the above frame localization dynamics for all agents we obtain a compact form
\begin{equation}\label{eq:finite-time_localization_compact_form}
\dot{\m{S}}(t)=-(\bar{\m{L}}\otimes \m{I}_4)\m{S}(t).
\end{equation}
where the matrix $\m{\bar{L}}=[\bar{l}_{ij}]\in \mb{R}^{n\times n}$ is defined as
$$
\bar{l}_{ij}=
\scriptstyle\begin{cases}
{0},~ (i,j)\in\mc{E},~i\neq j, ~\m{S}_i = \m{S}_j \text{ or }(i,j)\not\in\mc{E},~i\neq j \\
-{1}/{\abs{\m{S}_i-\m{S}_j}_F^\alpha},\qquad(i,j)\in\mc{E},~i\neq j,~\m{S}_i \neq \m{S}_j\\
\sum_{k\in\mc{N}_i}\bar{l}_{ik},\qquad i=j,i\in \mc{V}, 
\end{cases}
$$
which is a weighted Laplacian for the graph $\mc{G}$. 

Assume that $\mc{G}$ is undirected and connected. Then, $(\m{1}_n\otimes \m{I}_4)^T\dot{\m{S}}(t)=-(\m{1}_n\otimes \m{I}_d)^T(\m{\bar{L}}\otimes \m{I}_4)\m{S}(t)=\m{0}$. Thus, $(\m{1}_n\otimes \m{I}_4)^T\m{S}(t)$ is invariant under \eqref{eq:finite-time_localization_compact_form}. Let $\m{S}_c:=(1/n)(\m{1}_n\otimes \m{I}_4)^\top\m{S}(t)\in \mb{R}^{4\times 4}$, $\m{S}_i(t)=\m{S}_c+\boldsymbol{\delta}_i(t)$, and let $\boldsymbol{\delta}(t):=[\boldsymbol{\delta}_1^T,\ldots,\boldsymbol{\delta}_n^T]^T\in \mb{R}^{4n\times 4}$. Since $\m{S}_c$ is time-invariant, it follows that $\dot{\boldsymbol{\delta}}_i(t)=\dot{\m{S}}_i(t)$. Note that $\boldsymbol{\delta}_i-\boldsymbol{\delta}_j=\m{S}_i-\m{S}_j$.
\begin{Theorem}\label{thm:finiteTimeConverg}
Under the estimation law \eqref{eq:finite-time_law} and assume that $\mc{G}_o$ is a connected undirected graph, $\m{S}(t)$ globally asymptotically converges to $(\m{1}_n\otimes \m{I}_4)\text{Ave}\{\m{S}_i(0)\}_{i\in \mc{V}}$ in a finite time with settling time $T_c>0$ bounded by
\begin{equation*}
T_c\leq \frac{V(0)^{\alpha/2}}{\kappa\alpha}
\end{equation*} 
where $V(0)=\frac{1}{2}\sum_1^n||\boldsymbol{\delta_i}(0)||_F^2$, and $\kappa=(2\lambda_{2})^{\frac{2-\alpha}{2}}$ with $\lambda_{2}$ being the smallest nonzero eigenvalue of the Laplacian $\m{L}(\mc{G})$ associated with $\mc{G}$.
\end{Theorem}
\begin{pf}
Consider a Lyapunov candidate function
\begin{equation}\label{eq:LyapFuncUndirect}
V(t)=({1}/{2})\sum_{i=1}^{n}{\abs{\boldsymbol{\delta}_i}_F^2}=({1}/{2})\sum_{i=1}^{n}\text{tr}\big(\boldsymbol{\delta}^\top_i\boldsymbol{\delta}_i \big).
\end{equation}
Note that $V$ is radially unbounded, positive definite, continuously differentiable, and $V=0$ in $\mc{S}_o:=\big\{ \{\m{S}_i \}_{i\in\mc{V}}~|~\m{S}_i=\m{S}^c,~\forall i\in \mc{V}\big\}$.
The time derivative of $V$ along the trajectory of \eqref{eq:finite-time_localization_compact_form} is given as
\begin{align*}
\dot{V}(t)&=\textstyle\sum_{i=1}^{n}\text{tr}\big(\boldsymbol{\delta}^\top_i\boldsymbol{\dot{\delta}}_i \big)=\text{tr}\Big\{\textstyle\sum_{i=1}^{n}\boldsymbol{\delta}^\top_i\boldsymbol{\dot{\delta}}_i \Big\}
\\
&= -\text{tr}\Big\{\sum_{i=1}^{n}\boldsymbol{\delta}^\top_i\sum_{j\in\mc{N}_i}{\left(\boldsymbol{\delta}_i-\boldsymbol{\delta}_j\right)}/{\abs{\boldsymbol{\delta}_i-\boldsymbol{\delta}_j}_F^\alpha}\Big\}
\\
&= -\text{tr}\Big\{{\sum_{(i,j)\in\mc{E}}}{{\left(\boldsymbol{\delta}^\top_i-\boldsymbol{\delta}^\top_j\right)\left(\boldsymbol{\delta}_i-\boldsymbol{\delta}_j\right)}/{\abs{\boldsymbol{\delta}_i-\boldsymbol{\delta}_j}_F^\alpha}}\Big\}
\\
&= -\sum_{(i,j)\in\mc{E}}{{\abs{\boldsymbol{\delta}_i-\boldsymbol{\delta}_j}_F^2}/{\abs{\boldsymbol{\delta}_i-\boldsymbol{\delta}_j}_F^\alpha}}
\\
&=-\sum_{(i,j)\in\mc{E}}{\left(\abs{\boldsymbol{\delta}_i-\boldsymbol{\delta}_j}_F^2\right)^{(2-\alpha)/{2}}}
\\
&\leq -\Big(\sum_{(i,j)\in\mc{E}}\abs{\boldsymbol{\delta}_i-\boldsymbol{\delta}_j}_F^2\Big)^{(2-\alpha)/{2}} \numberthis \label{eq:dV_inequality_delta}
\\
&\leq -\Big[\text{tr}\big\{\sum_{(i,j)\in\mc{E}}(\boldsymbol{\delta}^\top_i-\boldsymbol{\delta}^\top_j)(\boldsymbol{\delta}_i-\boldsymbol{\delta}_j)\big\}\Big]^{(2-\alpha)/{2}}
\\
&\leq -\Big[\text{tr}\big\{\boldsymbol{\delta}^\top(\m{L}\otimes \m{I}_4)\boldsymbol{\delta}\big\}\Big]^{(2-\alpha)/{2}}
\\
&\leq -\Big[\sum_{k=1}^{4}[\boldsymbol{\delta}]^\top_k(\m{L}\otimes \m{I}_4)[\boldsymbol{\delta}]_k\Big]^{(2-\alpha)/{2}}\numberthis \label{eq:dV_inequality_delta_k}
\end{align*}
where $[\boldsymbol{\delta}]_k\in \mb{R}^{4n}$ denotes the $k$-th(1$\leq k \leq 4$) column vector of the matrix $\boldsymbol{\delta}\in \mb{R}^{4n\times 4}$ and the inequality \eqref{eq:dV_inequality_delta} is derived by applying Lemma \ref{lm:inequality} with $\frac{1}{2}<\frac{2-\alpha}{2}<1$. Under the assumption that $\mc{G}$ is connected undirected, $\m{L}\otimes \m{I}_d$ has $4$ zero eigenvalues and $\text{null}(\m{L}\otimes \m{I}_d)=\text{span}\{\text{Range}(\m{1}_n\otimes \m{I}_4)\}$. Since $[\boldsymbol{\delta}]_k\perp\text{Range}(\m{1}_n\otimes \m{I}_4)$, $\forall k = 1,\ldots,4$, one has
\begin{equation*}
[\boldsymbol{\delta}]^\top_k(\m{L}\otimes \m{I}_4)[\boldsymbol{\delta}]_k\geq \lambda_{2}[\boldsymbol{\delta}]^\top_k[\boldsymbol{\delta}]_k, ~\forall k=1,\ldots,4
\end{equation*} 
where $\lambda_{2}$ is the smallest nonzero eigenvalue of $\m{L}(\mc{G})$. Substituting the above inequalities into \eqref{eq:dV_inequality_delta_k} yields  
\begin{align*}
\dot{V}(t) &\leq -\Big[\lambda_{2}\sum_{k=1}^{4}\abs{[\boldsymbol{\delta}]_k}^2\Big]^{(2-\alpha)/{2}}
\\
&\leq -\lambda_{2}^{(2-\alpha)/{2}}\Big[\sum_{k=1}^{4}\sum_{i=1}^{n}\abs{[\boldsymbol{\delta}_i]_k}^2\Big]^{(2-\alpha)/{2}}
\\
&\leq -\lambda_{2}^{(2-\alpha)/{2}}\Big[\sum_{i=1}^{n}\text{tr}\big(\boldsymbol{\delta}^\top_i\boldsymbol{\delta}_i \big)\Big]^{(2-\alpha)/{2}}
\\
&\leq -\lambda_{2}^{(2-\alpha)/{2}}\big(2V(t)\big)^{(2-\alpha)/{2}}\\
&\leq -\kappa V(t)^{(2-\alpha)/{2}}, \numberthis \label{eq:dV_undirect_finite-time}
\end{align*}
where $\kappa=(2\lambda_{2})^{(2-\alpha)/{2}}$. It follows from Lemma \ref{lm:finiteTimeConverg} and \eqref{eq:dV_undirect_finite-time} that $V(t)$ converges to $0$ in finite time. In other words, $\{\m{S}_i\}_{i\in \mc{V}}$ converges to the invariant set $\mc{S}_o$. As a result, it follows that $\m{S}_i(t),~\forall i \in \mc{V}$, globally converges to $\m{S}^c=\text{Ave}\{\m{S}_i(0)\}_{i\in \mc{V}}$ with settling time $T_c\leq V(0)^{1-\frac{2-\alpha}{2}}/(\kappa(1-\frac{2-\alpha}{2}))=2V(0)^{\frac{\alpha}{2}}/(\kappa\alpha)$. \eprf
\end{pf}
\begin{Corollary}\label{corol:finite_time_global_convergence}
Assume that Assumption \ref{ass:measurements} holds and the interaction graph $\mc{G}$ is undirected and connected. Under the frame localization law \eqref{eq:finite-time_law}, if $\hat{\m{R}}_i^\top(t)$ is constructed from $\m{Q}_i(t)$ by the Gram-Schmidt procedure (GSOP) and $\hat{\m{p}}_i$ is computed by \eqref{eq:computated_pi_hat}, then there exist an unknown constant transformation  
$
\m{T}_c \in SE(3)
$
such that $\hat{\m{T}}_i(t)\rightarrow\m{T}_c^{-1}\m{T}_i(t)$ in a finite time, for all $i\in \mc{V}$, for almost all initializations of $\{\m{P}_i(0)\}_{i\in\mc{V}}$.
\end{Corollary}
\begin{pf}
The proof follows from the fact that $\m{S}_i(t)\rightarrow (\m{1}_n^\top\otimes \m{I}_4)\m{S}(0)=\text{Ave}\{\m{S}_i(0)\}_{i\in \mc{V}},\forall i \in \mc{V},$ globally asymptotically as $t\rightarrow T_c$ (Theorem \ref{thm:finiteTimeConverg}) and by following same lines as in the proof of Corollary \ref{corol:almost_global_convergence}. The estimated frame transformations of the agents, $\hat{\m{T}}_i$, converge almost globally asymptotically to the actual frame transformations in a finite-time, up to a common constant transformation $\m{T}_c(\m{R}_c,\m{q}_c)$, if the column vectors of $\m{Z}_0:=\left[\big(\m{R}_1(0)\m{Q}_1(0)\big)^\top,\ldots,\big(\m{R}_n(0)\m{Q}_n(0)\big)^\top\right]^\top\in \mb{R}^{3n\times 3}$ are not orthogonal to $\text{Range}(\m{1}_n\otimes \m{I}_3)$. \eprf
\end{pf}
\subsection{Relation between Actual Frame Transformations and Estimated Frame Transformations}
\begin{figure}[t]
\centering
\begin{tikzpicture}
\node (py) at (0,1.,0) [label=below left:$\Sigma$]{};
\node (pz) at (0,0,1.5) [label=left:$$]{};
\node (px) at (1.,-0.2) [label=below:$$]{};
\node[place,scale=0.6] (pi) at (1.7,1.5){};
\node[place,scale=0.6] (pj) at (3.2,.3){};
\node[place,scale=0.01] (pi_y) at (1.2,2.4)[label=right:$^i\Sigma$]{};

\draw[{line width=1pt}] (0,0) [frame] -- (px);
\draw[{line width=1pt}] (0,0) [frame] -- (py);
\draw[{line width=1pt}] (0,0) [frame] -- (pz);

\draw[{line width=.5pt}] (pi) [frame] -- (2.5,1.94);
\draw[{line width=.5pt}] (pi) [frame] -- (pi_y);
\draw[{line width=0.5pt}] (pi) [frame] -- (.7,1.6);

\draw[{line width=1pt}, blue] (0,0) [frame] -- (pi)  node [pos =0.6, yshift=-0.5ex, below] {$\m{T}_i$};
\draw[{line width=1pt}, blue] (0,0) [frame] -- (pj)  node [pos =0.5, yshift=-0ex, below] {$\m{T}_j$};
\draw[{line width=1pt}, blue] (pi) [frame] -- (pj)  node [pos =0.6, yshift=-0ex, above] {$\m{T}_{ij}$};

\draw[{line width=.5pt}] (pj) [frame] -- (3.2+0.1,.3+0.8) ;
\draw[{line width=0.5pt}] (pj) [frame] -- (3.2+1,0.1) node [above] {$^j\Sigma$};
\draw[{line width=0.5pt}] (pj) [frame] -- (3.2-0.6,-0.1);

\draw[{line width=.5pt},dashed,red] (pi) [frame] -- (2.7,3) node [pos =0.4, yshift=1.5ex, above] {\textcolor{black}{$-\m{q}_c$}};
\node[place,scale=0.6] () at (2.7,3){};
\draw[{line width=.5pt},dashed,red] (pj) [frame] -- (4.2,1.8) node [pos =0.4, yshift=2ex, above] {\textcolor{black}{$-\m{q}_c$}};
\node[place,scale=0.6] () at (4.2,1.8){};

\draw[red,dashed] (2.7,3) arc (40:60:4) node [pos =0.7, yshift=-0.5ex, below] {\textcolor{black}{$\m{R}_c^\top$}};
\draw[red,dashed] (4.2,1.8) arc (20:40:4)node [pos =0.3, yshift=1.5ex, above] {\textcolor{black}{$\m{R}_c^\top$}};

\draw[{line width=0.5pt}, blue, dashed] (0,0) [frame] -- (1.635,3.89)  node [pos =0.7, yshift=0ex, above left] {$\hat{\m{T}}_i$};
\node[place,scale=0.6] () at (1.635,3.89){};
\draw[{line width=0.5pt}, blue, dashed] (0,0) [frame] -- (3.5,3)  node [pos =0.9, yshift=-1ex, below] {$\hat{\m{T}}_j$};
\node[place,scale=0.6] () at (3.5,3){};

\draw[{line width=1pt}, blue] (1.635,3.89) [frame] -- (3.5,3)  node [pos =0.6, yshift=-0ex, above] {$\m{T}_{ij}$};
\end{tikzpicture}
\caption{Illustration of the relation between the estimated and actual frame transformations of the agents in $\mb{R}^3$.}
\label{fig:frame_relation}
\end{figure}
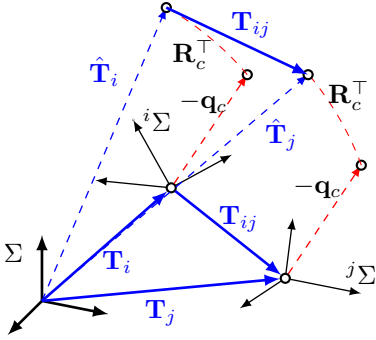

The relation between the estimates of the body motions, $\hat{\m{T}}_i$, and the actual body motions of the agents, $\m{T}_i$ is illustrated in Fig.  \ref{fig:frame_relation}. In particular, $\hat{\m{T}}_i,\forall i \in \mc{V},$ are obtained from $\m{T}_i$ by first translating the agents by $-\m{q}_c\in \mb{R}^3$ followed by a rotation $\m{R}_c^\top\in SO(3)$ about the origin of the global coordinate frame $\Sigma$. Note that the frame transformation $\m{T}_c(\m{R}_c,\m{q}_c)$ is an unknown constant and $\m{R}_{c}=\text{GSOP}(\m{Q}_c)$ where $\m{Q}_c$ and $\m{q}_c$ are computed as \eqref{eq:S_c}.
\section{Simulation}\label{sec:simulation}
\begin{figure*}
\centering
\begin{subfigure}[b]{0.2\textwidth}
\centering
\begin{tikzpicture}[scale=1.8]
\node[place] (2) at (1,0.) [label=right:$2$] {};
\node[place] (1) at (0,0) [label=left:$1$] {};
\node[place] (4) at (0,1.) [label=left:$4$] {};
\node[place] (3) at (1.,1.) [label=right:$3$] {};

\draw (2) [line width=1pt, -stealth'] -- node [left] {} (1);
\draw (3) [line width=1pt, -stealth'] -- node [left] {} (2);
\draw (4) [line width=1pt, -stealth'] -- node [left] {} (3);
\draw (2) [line width=1pt, -stealth'] -- node [left] {} (4);
\end{tikzpicture}
\caption{$\mc{G}$.}
\label{fig:sim_graph}
\end{subfigure}
\begin{subfigure}[b]{0.38\textwidth}
\centering
\includegraphics[height=4.8cm]{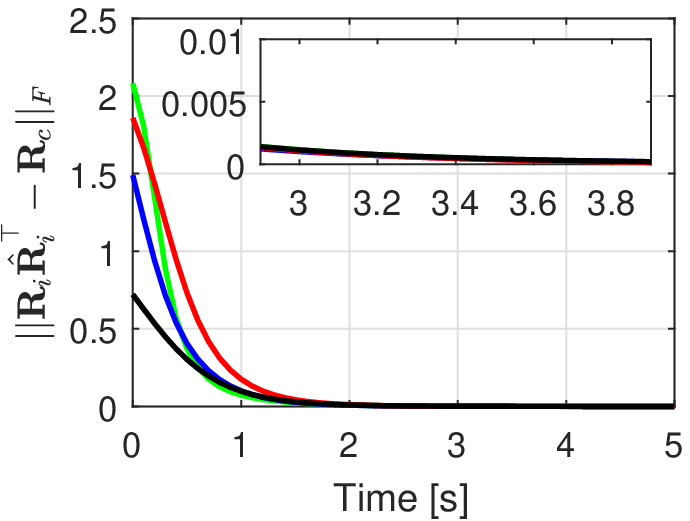}
\caption{Orientation errors.}
\label{fig:orient_err}
\end{subfigure}
\begin{subfigure}[b]{0.38\textwidth}
\centering
\includegraphics[height=4.8cm]{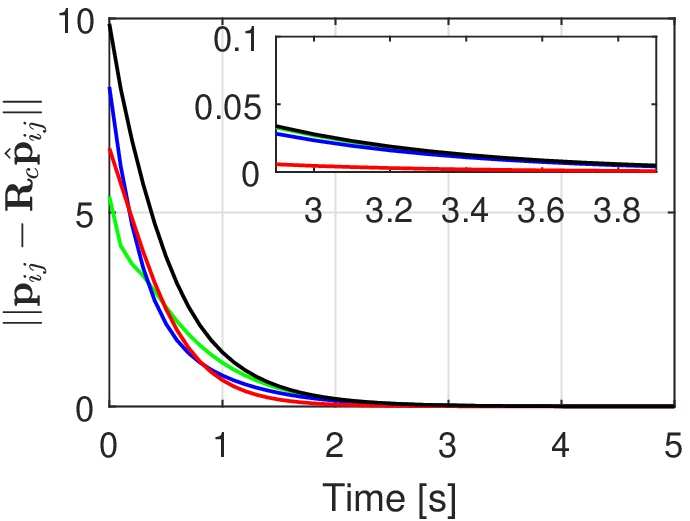}
\caption{Position errors.}
\label{fig:pos_err}
\end{subfigure}
\caption{Frame localization induced norm errors of four agents under localization law \eqref{eq:asymptotic_law}.}
\label{fig:sim_frame_local}
\end{figure*}
\begin{figure*}
\centering
\begin{subfigure}[b]{0.2\textwidth}
\centering
\begin{tikzpicture}[scale=1.8]
\node[place] (2) at (1,0.) [label=right:$2$] {};
\node[place] (1) at (0,0) [label=left:$1$] {};
\node[place] (4) at (0,1.) [label=left:$4$] {};
\node[place] (3) at (1.,1.) [label=right:$3$] {};

\draw (2) [line width=1pt] -- node [left] {} (1);
\draw (3) [line width=1pt] -- node [left] {} (2);
\draw (4) [line width=1pt] -- node [left] {} (3);
\draw (1) [line width=1pt] -- node [left] {} (4);
\end{tikzpicture}
\caption{$\mc{G}$.}
\label{fig:sim_graph_finite}
\end{subfigure}
\begin{subfigure}[b]{0.38\textwidth}
\centering
\includegraphics[height=4.8cm]{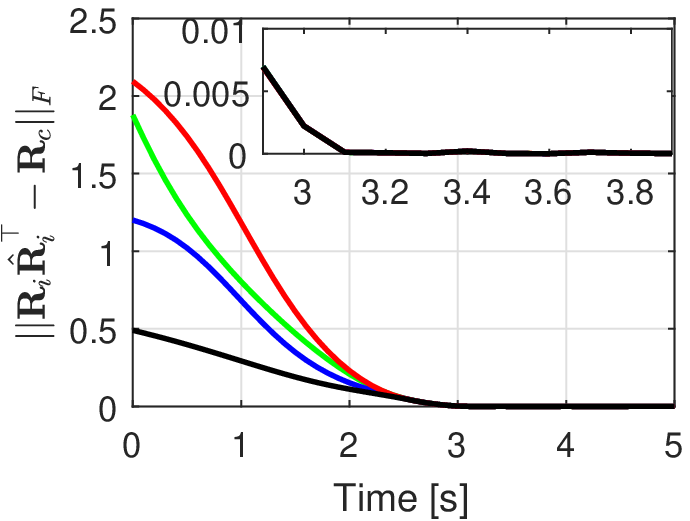}
\caption{Orientation errors.}
\label{fig:orient_err_finite}
\end{subfigure}
\begin{subfigure}[b]{0.38\textwidth}
\centering
\includegraphics[height=4.8cm]{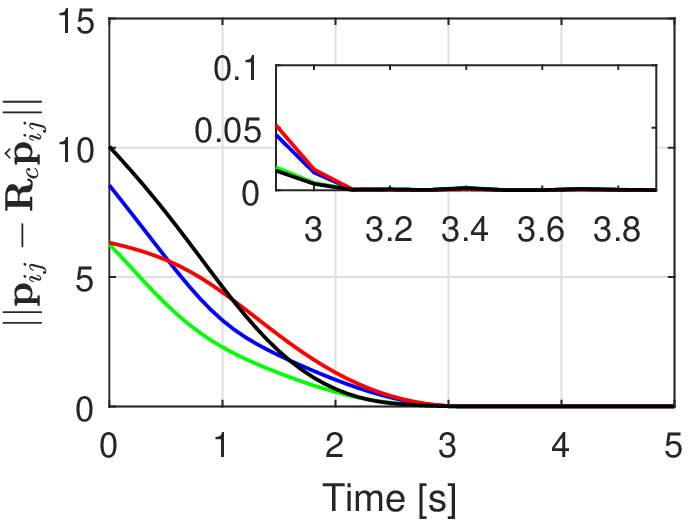}
\caption{Position errors.}
\label{fig:pos_err_finite}
\end{subfigure}
\caption{Frame localization induced norm errors of four agents under finite-time localization law \eqref{eq:finite-time_law}.}
\label{fig:sim_frame_local_finite}
\end{figure*}

Consider a system of four mobile agents in $\mb{R}^3$. The initial orientations of the agents are chosen randomly and their initial positions are given as: $\m{p}_1(0)=[0,0,0]^\top,~ \m{p}_2(0)=[4,0,0]^\top,~ \m{p}_3(0)=[0,0,4]^\top,$ and $\m{p}_4(0)=[0,4,0]^\top$. 
The agents travel with angular velocities: $\boldsymbol{\omega}_1^1=[0.3,0,0]^\top,~\boldsymbol{\omega}_2^2=[0,0.3,0]^\top,~\boldsymbol{\omega}_3^3=[0,0,0.3]^\top,~\boldsymbol{\omega}_4^4=\m{0}$, and linear velocities: $\m{v}_1=[1,0,0]^\top,~\m{v}_2=[0,1,0]^\top,~\m{v}_3=[0,0,1]^\top$ and $\m{v}_4=[1,1,0]^\top$, respectively. 

We provide simulation results for frame localization of the system of four agents under the localization law \eqref{eq:asymptotic_law} and the finite-time localization law \eqref{eq:finite-time_law} with the corresponding interaction graph topologies in Figs. \ref{fig:sim_frame_local} and \ref{fig:sim_frame_local_finite}, respectively. In both cases, it is observed that the estimated poses of the agents asymptotically converge to the actual poses as $t\rightarrow\infty$ since the induced norms of the orientation errors, $||\m{R}_i\hat{\m{R}}_i^\top-\m{R}_c||_F\rightarrow 0$ (See Figs. \ref{fig:orient_err} and \ref{fig:orient_err_finite}), and the norms of the position errors $||(\m{p}_j-\m{p}_i)-\m{R}_c(\hat{\m{p}}_j-\hat{\m{p}}_i)||\rightarrow0$ (See Figs. \ref{fig:pos_err} and \ref{fig:pos_err_finite}). 

To show the advantage of the finite-time localization scheme \eqref{eq:finite-time_law} over the localization protocol with asymptotic stability property \eqref{eq:asymptotic_law} the zoomed plots of the orientation and position estimation errors in the time interval $t\in [3,4]$ are included. It is shown that the finite-time frame localization scheme \eqref{eq:finite-time_law} has fast convergence time (about $3.2$ seconds), whereas, the convergence time of the estimated poses under \eqref{eq:asymptotic_law} is infinitely long.
\section{Conclusion}\label{sec:Conclusion}

In this work, we presented two frame localization schemes for estimating the trajectories of rigid-body motions of multi-agent systems in $\mb{R}^3$. Under the first localization law, the estimated frame transformations of the agents converge to the actual frame transformations almost globally and exponentially up to an unknown constant transformation bias. Whereas, under the second frame localization protocol, the estimated frame transformations of the agents converge to the actual frame transformations in a finite time up to an unknown constant transformation bias. Simulation results were provided to support the theoretical analysis. There are several possible applications of the proposed frame localization schemes in the network pose formation control and the multi-agent cooperative estimation and which are left as the future work.

\bibliography{ifacconf}             
                                                   







\appendix
\section{The Gram-Schmidt Orthonormalization Process (GSOP)}\label{app:gsop}
For a set of $d$ independent vectors $\mc{Z}=\{\m{z}_1,\m{z}_2,\ldots,\m{z}_d\}$ in $\mb{R}^d$, the GSOP, which constructs $d$ orthonormal column vectors of $\m{Q}=[\m{q}_1~\ldots\m{q}_d]\in\text{SO}(d)$ from $\mc{Z}$, is defined as follows,
\begin{align*}
&& \m{v}_1&:=\m{z}_1, && \shoveright{\m{q}_1:=\m{v}_1/\abs{\m{v}_1},}\\
&& \m{v}_2&:=\m{z}_2-\langle\m{z}_2,\m{q}_1 \rangle\m{q}_1, && \shoveright{\m{q}_2:=\m{v}_2/{\abs{\m{v}_2}},}\\
&& & ~~~\ldots && ~~~~~\ldots\\
&& \m{v}_d&:=\m{z}_d-\textstyle\sum\nolimits_{k=1}^{d-1}\langle\m{z}_d,\m{q}_k \rangle\m{q}_k, && \shoveright{\m{q}_d:=\alpha{\m{v}_d}/{\abs{\m{v}_d}},}
\end{align*}
where $\langle\cdot,\cdot\rangle$ denotes the inner product, and the coefficient $\alpha$ is chosen such that $\text{det}(\m{Q})=1$ as
$\alpha:=\text{sign}\left( \text{det}\left( [\m{q}_1, \ldots, \m{q}_{d-1},{\m{v}_d}/{\abs{\m{v}_d}}] \right)  \right).$ 
If the set $\mc{Z}$ contains a linearly dependent vector (which linearly depends on one or more vectors in $\mc{Z}$), there exists $\m{v}_i=\m{0},~i\in \{1,\ldots,d\}$, and hence $\m{q}_i=\m{0}$.
\begin{Lemma}\cite[Lemma 2]{Quoc2018cdc}\label{lm:invariantOrthogonalization}
Consider the rotation transformation $\m{Q}_i=\m{R}_i^\top\m{Q}_{Si}, \m{R}_i\in SO(3), \m{Q}_i$ and $\m{Q}_{Si}\in \mb{R}^{3\times 3}$. If $\m{\hat{R}}_i^\top$ and $\m{R}_{Si}$ are derived from $\m{Q}_{i}$ and $\m{Q}_{Si}$ by the above GSOP, respectively, then there holds
\begin{equation}\label{eq:estimatedOrient}
\m{\hat{R}}_i^\top=\m{R}_i^\top\m{R}_{Si}.
\end{equation}
\end{Lemma}
\section{Finite-Time Convergence Theory}\label{app:finite_time}

\begin{Lemma}\cite{Hardy1952}\label{lm:inequality}
If $\xi_1,\ldots,\xi_d\geq 0$ and $0\leq p \leq 1$, then
\begin{equation*}
\left( \textstyle\sum_{i=1}^{d}\xi_i\right)^p\leq \textstyle\sum_{i=1}^{d}\xi_i^p.
\end{equation*}
\end{Lemma}

A condition for finite-time convergence of continuous-time systems is given by the following lemma.
\begin{Lemma}(\cite{Bhat2000}).\label{lm:finiteTimeConverg}
Suppose that there exists a positive-definite and continuous function $V(x,t):\mb{R}^n\times [0,\infty)\rightarrow \mb{R}$. If there exists $\kappa>0,~\alpha\in (0,1)$, and open neighborhood $\mc{U}_0\in \mb{R}^d$ of the origin such that
\begin{equation*}
\dot{V}+\kappa V^\alpha\leq 0,~\forall \m{x}\in\mc{U}_0\setminus \{\m{0}\},  
\end{equation*}
then $V=0$ for $t\geq T$, with the settling time $T\leq V^{1-\alpha}(0)/(\kappa(1-\alpha))$.
\end{Lemma}
\end{document}